\documentclass{l4dc2024}

\title[Model-Free NN Verification]
{Model-Free Verification for Neural Network Controlled Systems
}
\usepackage{times}
\usepackage{mathrsfs}
\newtheorem{problem}{Problem}
\newtheorem{assumption}{Assumption}
\author{%
 \Name{Han Wang} \Email{han.wang@eng.ox.ac.uk}\\
 \Name{Zuxun Xiong} \Email{zuxun.xiong@eng.ox.ac.uk}\\
 \Name{Liqun Zhao} \Email{liqun.zhao@eng.ox.ac.uk}\\
 \Name{Antonis Papachristodoulou} \Email{antonis@eng.ox.ac.uk}\\
 \addr Department of Engineering Science, University of Oxford, Parks Road, OX1 3PJ, Oxford, UK%
}

\begin{document}

\maketitle

\begin{abstract}%
Neural network controllers have shown potential in achieving superior performance in feedback control systems. Although a neural network can be trained efficiently using deep and reinforcement learning methods, providing formal guarantees for the closed-loop properties is challenging. The main difficulty comes from the nonlinear activation functions. One popular method is to use sector bounds on the activation functions resulting in a robust analysis. These methods work well under the assumption that the system dynamics are perfectly known, which is, however, impossible in practice. In this paper, we propose data-driven semi-definite programs to formally verify stability and safety for a neural network controlled linear system with unknown dynamics. The proposed method performs verification directly from end-to-end without identifying the dynamics. Through a numerical example, we validate the efficacy of our method on linear systems with controller trained by imitation learning.
\end{abstract}


\section{Introduction}
Neural Networks (NNs) have attracted an increasing attention in both academia and industry over the past decade. Following advances in computational hardware, data collection and storage techniques, and deep learning, NNs with a large amount of parameters can be trained for accomplishing complex tasks. They have shown great power in a range of applications such as face recognition \citep{parkhi2015deep}, natural language processing \citep{chowdhary2020natural}, weather prediction \citep{ren2021deep}, etc. One of the most exciting applications of NNs is in control theory, where autonomous devices and processes are controlled by neural networks for advanced behaviours.

In control theory, rigorous linear system theories have been developed for stability \citep{kailath1980linear}, robustness \citep{zhou1998essentials}, and optimality \citep{khalil1996robust}. For nonlinear systems, several convex optimization based methods have been proposed such as optimal control \citep{lewis2012optimal}, model predictive control \citep{kouvaritakis2016model}, and state-feedback linear matrix inequalities \citep{boyd1994linear}. However, all these optimization based methods are hard to be applied for systems of high dimension and complex nonlinearities. From the optimization point of view, the optimization problem becomes large and potentially nonconvex. Moreover, when controlling a real autonomous system in the real-world, control should be designed to certify more properties, such as safety \citep[Chapter 4.B]{annaswamy2023control}. More recently, NN controllers designed by imitation learning (IL) and reinforcement learning \citep{zhao2023barrier} have been proposed and applied to nonlinear control applications, such as racing drones \citep{kaufmann2023champion}. However, before these NN controllers are implemented, verification of critical properties such as safety and stability should be performed. This is because learning methods are usually hard to provide guarantees. 

Performing an accurate verification for an NN controlled closed-loop system is challenging due to the nonlinear activation functions present in NNs. A class of methods to achieve this is to find bounds for the activation functions. Initial attempts focussed on finding the worst-case outcome of each layer \citep{gowal2018effectiveness}. Later on, tighter bounds using quadratic constraints based on sectors and slopes have been proposed \citep{fazlyab2020safety}. This method has shown to be efficient for various activation functions, such as ReLU, leaky ReLU, tanh, sigmoid, etc. Using multiple sectors to bound the activation functions has shown better verification performance in the reachability problem \citep{newton2021neural}. The obtained bound can be encoded as semi-definite constraints in the verification program \citep{yin2021imitation}. Exploiting the sparsity in the NN can accelerate the solutiuon of the SDP using chordal decomposition \citep{newton2023sparse}. However, all these methods are based on the assumption that the system dynamics are known. 

In this work, we consider the problem of verifying a neural network controlled linear system, when the system dynamics are unknown. The idea is to collect input and state data that are persistently exciting to represent the system dynamics using Willems's fundamental lemma \citep{willems2005note}, and formulate a data-driven sum-of-squares program for verification. Properties that are verified in this paper involve stability and safety.

In Section \ref{sec:prob-formulation} we formulate the verification problem for an unknown linear system, and provide the necessary preliminaries. The main results are presented in Section \ref{sec:nn-verification}. Section \ref{sec:simulation} shows simulation result, and Section \ref{sec:conclusion} concludes the paper.

\section{Problem Formulation}
\label{sec:prob-formulation}
Consider a feedback system consisting of an unknown plant $G$ and a feed-forward fully connected Neural Network controller $\pi:\mathbb{R}^{n_x}\to \mathbb{R}^{n_u}$, which satisfies $\pi(0)=0$. Even though the dynamics of $G$ are unknown, this is assumed to be a discrete-time linear system
\begin{equation}\label{eq:system}
    x(k+1) = A_Gx(k)+B_G\pi(k),
\end{equation}
where $x(k)\in\mathbb{R}^{n_x}$ denotes the state, $u(k)\in\mathbb{R}^{n_u}$ denotes the control input and $A_G\in\mathbb{R}^{n_x\times n_x}$ and $B_G\in\mathbb{R}^{n_x\times n_u}$ are unknown matrices which are assumed to be time-invariant. The neural network controller $\pi$ is defined as follows:
\begin{subequations}\label{eq:nn}
    \begin{align}
        &w^0(k) = x(k),\\
        &v^i(k) = W^iw^{i-1}(k)+b^i,~\mathrm{for}~i=1,\ldots,l,\\
        &w^i(k)=\phi^i(v^{i}(k)),i=1,\ldots,l,\\
        &\pi(k)=W^{l+1}w^l(k)+b^{l+1},
    \end{align}
\end{subequations}
where $W^i\in\mathbb{R}^{n_{i}\times n_{i-1}}$ is the weight matrix, $b^{i}\in\mathbb{R}^{n_i}$ is the bias vector of the $i^{\mathrm{th}}$ layer, $v^i$ and $w^i$ are the corresponding vectorized input and output, respectively. The input of the NN is $w^0(k)=x(k)$, which is the state of the system at time $k$. $\phi^i:\mathbb{R}^{i}\to \mathbb{R}^i$ is the vector of nonlinear activation functions in the $i^{\mathrm{th}}$ layer, defined as
\begin{equation}
    \phi^i(v^i) = [\varphi(v^i_1),\ldots,\varphi(v^i_{n_i})]^\top,
\end{equation}
where $\varphi:\mathbb{R}\to \mathbb{R}$ is the activation function, such as ReLU $\varphi(\nu) := \max(0,\nu)$, tanh $\varphi(\nu):=\tanh(\nu)$, sigmoid $\varphi(\nu) :=\frac{1}{1+e^{-\nu}}$, etc, for the NN. The verification problem is then formulated as follows.

\begin{problem}[Stability Verification]\label{prob:stability}
    Consider system \eqref{eq:system}, where $A_G$, $B_G$ are unknown, $\pi$ is a given neural network controller as in \eqref{eq:nn}. Verify stability of system \eqref{eq:system} around the equilibrium $x_*=0$.
\end{problem}

\begin{problem}[Finite-Time Safety Verification]\label{prob:safety}
    Consider system \eqref{eq:system}, where $A_G$, $B_G$ are unknown, $\pi$ is a given neural network controller as in \eqref{eq:nn}. Let $\mathcal{X}$ be the set of initial states, i.e. $x(0)\in\mathcal{X}$, $\mathcal{S}$ be the set of safe states. Verify that $x(k)\in\mathcal{S}$ for $k=1,\ldots,T$, and $\forall x(0)\in\mathcal{X}$.
\end{problem}

\begin{problem}[Invariance Verification]\label{prob:invariance}
Consider system \eqref{eq:system}, where $A_G$, $B_G$ are unknown, $\pi$ is a given neural network controller as in \eqref{eq:nn}. Let $\mathcal{B}\subseteq \mathbb{R}^{n_x}$ be a given set. Verify if $A_Gx+B_G\pi\in\mathcal{B}$ for any $x\in\mathcal{B}$.
    
\end{problem}

Before illustrating our method, we provide necessary preliminaries for analyzing the unknown plant $G$ and nonlinearities in the NN controller $\pi$.

\subsection{Data-Driven Open Loop Representation}
Consider a $T$-long time series of inputs, states, and successor states data of system \eqref{eq:system}:
\begin{subequations}
    \begin{align}
    &U_0:=\begin{bmatrix}
        u(0)&u(1)&\ldots&u(K-1)
    \end{bmatrix}\in\mathbb{R}^{n_u\times K}\\
    &X_0:=\begin{bmatrix}
        x(0)&x(1)&\ldots&x(K-1)
    \end{bmatrix}\in\mathbb{R}^{n_x\times K}\\
    &X_1:=\begin{bmatrix}
        x(1)&x(2)&\ldots&x(K)
    \end{bmatrix}\in\mathbb{R}^{n_x\times K}
\end{align}
\end{subequations}
We assume that the data are collected in an accurate way without noise corruptions. We also assume that the data is sufficiently rich, which indicates
\begin{equation}
    \mathrm{rank}\left(\begin{bmatrix}
        U_0\\X_0
    \end{bmatrix}\right)=n_u+n_x.
\end{equation}
Then, system \eqref{eq:system} has the following equivalent data-driven representation \citep[Theorem 1]{de2019formulas}:
\begin{equation}
    x(k+1)=X_1\begin{bmatrix}
        U_0\\X_0
    \end{bmatrix}^\dag
    \begin{bmatrix}
    \pi(k)\\
    x(k)
    \end{bmatrix},
\end{equation}
where $\dag$ represents the right inverse.
Following this, we immediately have
\begin{equation}
    \begin{bmatrix}
        B_G&A_G
    \end{bmatrix} = X_1\begin{bmatrix}
        U_0\\X_0
    \end{bmatrix}^\dag.
\end{equation}

\subsection{Neural Network Representation}
In order to deal with the nonlinearities in the NN controller \eqref{eq:nn}, nonlinear isolation techniques has been proposed \citep{yin2021stability}. The $l$-layer neural network controller \eqref{eq:nn} can be written as
\begin{subequations}\label{eq:nnsystem}
\begin{align}
    \begin{bmatrix}
        \pi(k)\\
        v_\phi(k)
    \end{bmatrix}&=N\begin{bmatrix}
        x(k)\\w_\phi(k)
    \end{bmatrix}\\
    w_\phi(k)&=\phi(v_\phi(k)),
\end{align}
\end{subequations}
where $v_\phi(k)$ and $w_\phi(k)$ denotes the stacked inputs and outputs of all the neurons in the neural network at time $k$, i.e.,
\begin{equation}
    v_\phi(k):=\begin{bmatrix}
        v^1(k)^\top&\ldots&v^l(k)
    \end{bmatrix}^\top \in\mathbb{R}^{n_\phi},w_\phi(k):=\begin{bmatrix}
        w^1(k)^\top&\ldots&w^l(k)
    \end{bmatrix}\in\mathbb{R}^{n_\phi},
\end{equation}
 $n_\phi=n_1+\ldots,n_l$. Similarly, the stacked activation function $\phi:\mathbb{R}^{n_\phi}\to \mathbb{R}^{n_\phi}$ is defined as
 \begin{equation}
     \phi(v_\phi):=\begin{bmatrix}
         \phi^1(v^1)^\top & \ldots & \phi^l(v^l)
     \end{bmatrix}^\top. 
 \end{equation}
The matrix $N$ can be derived from \eqref{eq:nn} as follows
\begin{equation}\label{eq:}
    N: = \left[ {\begin{array}{*{20}{c}}
0&\vline& 0&0& \ldots &{{W^{l + 1}}}\\
\hline
{{W^1}}&\vline& 0& \cdots &0&0\\
0&\vline& {{W^2}}& \cdots &0&0\\
 \vdots &\vline&  \vdots & \ddots & \vdots & \vdots \\
0&\vline& 0& \cdots &{{W^l}}&0
\end{array}} \right]:=\begin{bmatrix}
    N_{ux}&N_{uw}\\N_{vx}&N_{vw}
\end{bmatrix}.
\end{equation}

\subsection{Sector Constraints}
Sector constraints are commonly used to bound the nonlinear activation function $\phi$ with quadratic constraints. For one activation function $\varphi$, the sector bound is defined as follows.
\begin{definition}
    Given $\alpha\le \beta$. The function $\varphi:\mathbb{R}\to \mathbb{R}$ lies in the sector $[\alpha,\beta]$ if 
    \begin{equation}
        (\varphi(\nu)-\alpha(\nu))(\beta(\nu)-\varphi(\nu))\ge 0,\forall \nu\in\mathbb{R}.
    \end{equation}
\end{definition}

The activation functions tanh, ReLU, sigmoid, can all be bounded by sectors. The ReLU function $\varphi(\nu)=\max(0,\nu)$ and tanh $\varphi(\nu)=\tanh(\nu)$ are in sector $[0,1]$. From a geometric point of view, the function $\varphi(\nu)$ is in the sector $[\alpha,\beta]$ if it lies between the lines $\alpha \nu$ and $\beta\nu$. These lines intersect at the origin, which can be regarded as the centre of the corresponding sector constraint. As an alternative, offset sector constraints can be defined if the centre is shifted from $(0,\varphi(0))$ to a certain point $(\nu_*,\varphi(\nu_*))$.
\begin{definition}
    Given $\alpha\le \beta\in\mathbb{R}$, $\nu_*\in\mathbb{R}$. The function $\varphi:\mathbb{R}\to \mathbb{R}$ satisfies the offset sector $[\alpha,\beta]$, around the point $\nu_*,\varphi(\nu_*)$ if 
    \begin{equation}
        ((\varphi(\nu_*)-\varphi(\nu))-\alpha(\nu_*-\nu))(\beta(\nu_*-\nu)-(\varphi(\nu_*)-\varphi(\nu)))\ge 0,\forall \nu\in\mathbb{R}.
    \end{equation}
\end{definition}
For example, the sigmoid $\varphi(\nu)= \frac{1}{1+e^{-\nu}}$ has been shown to satisfy the offset sector with $[\nu_*,\varphi(\nu_*)]=[0,0.5]$, $[\alpha,\beta]=[0,0.25]$.

Let now $w_{\phi}^i=\varphi(v_{\phi}^i)$ denote the output of the $i^{\mathrm{th}}$ activation function in $\phi(v_\phi)$. For $i=1,\ldots,n_{\phi}$, assume $\varphi(v_{\phi}^i)$ satisfies the offset sector constraints around $[v_{*}^i,\varphi(v_*^i)]$, within sector $[\alpha_{\phi}^i,\beta_{\phi}^i]$. Define $\alpha_{\phi}=[\alpha_\phi^1,\ldots,\alpha_{\phi}^{n_\phi}]$, $\beta_{\phi}=[\beta_\phi^1,\ldots,\beta_{\phi}^{n_\phi}]$, $v_*=[v^1_*,\ldots,v^{n_\phi}_*]$, $w_*=[w^1_*,\ldots,w^{n_\phi}_*]$. The following lemma bounds the vectored activation function $\phi(v_{\phi})$.

\begin{lemma}[\protect {\citep[Lemma 1]{yin2021stability}}]\label{lem:sector}
    Let $\alpha_\phi$, $\beta_{\phi}\in\mathbb{R}^{n_\phi}$ be given with $\alpha_{\phi}\le \beta_{\phi}$, and $w_*:=\phi(v_*)$. Assume $\phi$ element-wisely satisfies the offset sector $[\alpha_{\phi},\beta_{\phi}]$ around the point $(v_*,w_*)$ for all $v_\phi\in\mathbb{R}^{n_\phi}$. Then, for any $\lambda\in\mathbb{R}^{n_\phi}$ with $\lambda\ge 0$, and for all $v_\phi\in\mathbb{R}^{n_\phi}$, $w_\phi=\phi(v_\phi)$, we have
    \begin{align}\label{eq:sector}
        \begin{bmatrix}
            v_\phi-v_*\\
            w_\phi-w_*
        \end{bmatrix}^\top
        \begin{bmatrix}
            -2A_\phi B_\phi\Lambda&(A_\phi+B_\phi)\Lambda\\
            (A_\phi+B_\phi)\Lambda&-2\Lambda
        \end{bmatrix}
        \begin{bmatrix}
            v_\phi-v_*\\
            w_\phi-w_*
        \end{bmatrix}\ge 0,
    \end{align}
    where $A_\phi$ = diag($\alpha_\phi$), $B_\phi$ = diag($\beta_\phi$), and $\Lambda$ = diag($\lambda$).
\end{lemma}

\section{Neural Network Verification}
\label{sec:nn-verification}

In this section, convex programs to solve the Problems \ref{prob:stability} and \ref{prob:safety} are provided. To verify the neural network controlled closed-loop system with an unknown plant, we assume that measurements for state signal $x(k)$ and input signal $\pi(k)$ are obtainable. Consider the open-loop system $\dot x =A_Gx+B_Gu$, let 
\begin{align}\label{eq:data}
    &U_0:=\begin{bmatrix}
        u(0)&u(1)&\ldots&u(K-1)
    \end{bmatrix}\in\mathbb{R}^{n_u\times K}\\
    &X_0:=\begin{bmatrix}
        x(0)&x(1)&\ldots&x(K-1)
    \end{bmatrix}\in\mathbb{R}^{n_x\times K}\\
    &X_1:=\begin{bmatrix}
        x(1)&x(2)&\ldots&x(K)
    \end{bmatrix}\in\mathbb{R}^{n_x\times K}
\end{align} be input-output data sampled from $K$ independent experiments. 

Assumptions are given as follows.
\begin{assumption}\label{ass:data-exciting}
\begin{equation}\label{eq:exciting}
    \mathrm{rank}\left(\begin{bmatrix}
        U_0\\X_0
    \end{bmatrix}\right)=n_x+n_u,\quad \mathrm{rank}\left(X_1\right)=n_x.
\end{equation}
\end{assumption}

\begin{assumption}\label{ass:safe}
    The set $\mathcal{X}\subset\mathbb{R}^{n_x}$ is a non-empty polytope, defined by 
    \begin{equation}\label{eq:input-set}
        \mathcal{X}:=\left\{x\in\mathbb{R}^{n_x}:\bigcap_{i=1}^oc_i^\top (x-x^{ini})+1\ge 0\right\},
    \end{equation}
    where $c_i\in\mathbb{R}^{n_x+1}$ for $i=1,\ldots,o$, $x^{ini}\in\mathcal{X}\subset \mathbb{R}^{n_x}$ is a point inside $\mathcal{X}$. The set $\mathcal{S}\subset \mathbb{R}^{n_x}$ is also a non-empty polytope, defined by
    \begin{equation}\label{eq:safe-set}
        \mathcal{S}:=\left\{x\in\mathbb{R}^{n_x}:\bigcap_{i=1}^pd_i^\top (x-x^{safe})+1\ge 0\right\},
    \end{equation}
   where $d_i\in\mathbb{R}^{n_x}$ for $i=1,\ldots,p$, $x^{safe}\in\mathcal{X}\subset \mathbb{R}^{n_x}$ is a point inside $\mathcal{S}$. 
\end{assumption}

\begin{assumption}\label{ass:invariance}
The set $\mathcal{B}\subset\mathbb{R}^{n_x}$ is a non-empty polytope, defined by
\begin{equation}\label{eq:invariant}
    \mathcal{B}:=\left\{x\in\mathbb{R}^{n_x}:\bigcap_{i=1}^q e_i^\top (x-x^{inv})+1 \ge 0\right\},
\end{equation}
where $e_i\in\mathbb{R}^{n_x}$ for $i=1,\ldots,q$, $x^{inv}\in\mathcal{B}\subset \mathbb{R}^{n_x}$ is a point inside $\mathcal{B}$.
    
\end{assumption}
We can now present the main result of this paper.

\subsection{Stability Verification}
To obtain convex conditions for stability verification, we perform a loop transformation for the NN controller \eqref{eq:nn} following \citep[Section IV.A]{yin2021imitation}. The equivalent NN representation can be described by
\begin{subequations}\label{eq:NN-loop-transform}
    \begin{align}
        \begin{bmatrix}
            \pi(k)\\
            v_\phi(k)
        \end{bmatrix}&=
        \tilde N \begin{bmatrix}
            x(k)\\
            z_\phi(k)
        \end{bmatrix}\\
        z_\phi(k)&=\tilde \phi(v_\phi(k)).
    \end{align}
\end{subequations}
The transformed nonlinearity $\tilde \phi$ has been normalized, thus satisfies
\begin{equation}
    \begin{bmatrix}
        v_\phi\\z_\phi
    \end{bmatrix}^\top 
    \begin{bmatrix}
        \Lambda&0\\0&-\Lambda
    \end{bmatrix}
    \begin{bmatrix}
        v_\phi\\z_\phi
    \end{bmatrix}\ge 0,\quad \forall v_\phi\in[\underline{v},\overline{v}].
\end{equation}
The transformed matrix $\tilde N$ is derived by
\begin{equation}
    \tilde N=\begin{bmatrix}
        N_{ux}+C_2(I-C_4)^{-1}N_{vx}&C_1+C_2(I-C_4)^{-1}C_3\\(I-C_4)^{-1}N_{vx}&(I-C_4)^{-1}C_3
    \end{bmatrix}:=
    \begin{bmatrix}
        \tilde N_{ux}&\tilde N_{uz}\\
        \tilde N_{vx}&\tilde N_{vz}
    \end{bmatrix}
\end{equation}
where
\begin{equation}
    C_1=N_{uw}\frac{B_\phi-A_\phi}{2},C_2=N_{uw}\frac{A_\phi+B_\phi}{2},C_3=N_{vw}\frac{B_\phi-A_\phi}{2},C_4=N_{vw}\frac{A_\phi+B_\phi}{2}.
\end{equation}
\begin{theorem}\label{th:stability}
    Consider System \eqref{eq:system} with an unknown plant $G$ and an NN controller $\pi$ as in \eqref{eq:nn}. Let Assumption \ref{ass:data-exciting} hold. Let $\alpha_\phi$, $\beta_{\phi}$, $v_*\in\mathbb{R}^{n_\phi}$ be given with $\alpha_{\phi}\le \beta_{\phi}$, and $w_*:=\phi(v_*)$. Assume $\phi$ element-wisely satisfies the offset sector $[\alpha_\phi,\beta_\phi]$, centered on the point $(v_*,w_*)$. If there exist $L_1\in\mathbb{R}^{K\times n_x}$, $L_2\in\mathbb{R}^{K\times n_\phi}$, $Q_1\in\mathbb{S}^{n_x}_{++}$, $Q_2\in\mathbb{S}^{n_\phi}_{++}$, $Q_2$ is a diagonal matrix, such that
    \begin{subequations}\label{eq:data-stability-verification}
    \begin{equation}\label{eq:stability-LMI}
            H:=\begin{bmatrix}
                Q_1&0&L_1^\top X_1^\top&Q_1\tilde{N}_{vx}^\top\\
                0&Q_2&L_2^\top X_1^\top&Q_2\tilde{N}_{vz}^\top\\
                X_1L_1&X_1L_2&Q_1&0\\
                \tilde{N}_{vx}Q_1&\tilde{N}_{vz}Q_2&0&Q_2
            \end{bmatrix}\succ 0
    \end{equation}
    \begin{equation}\label{eq:data-representation-matrix}
        \begin{bmatrix}
            \tilde{N}_{ux}Q_1\\Q_1
        \end{bmatrix}=\begin{bmatrix}
            U_0\\X_0
        \end{bmatrix}L_1,\begin{bmatrix}
            \tilde{N}_{uz}Q_2\\0
        \end{bmatrix}=\begin{bmatrix}
            U_0\\X_0
        \end{bmatrix}L_2
    \end{equation}\label{eq:data-stability-matrix}
        \end{subequations}
then the Closed-loop System \eqref{eq:system} is asymptotically stable around $x_*=0$.
\end{theorem}
\begin{proof}
Using \citep[Lemma 2]{yin2021imitation}, the NN controlled system \eqref{eq:system} is asymptotically stable around $x_*$ if there exists a matrix $P\in\mathbb{S}^{n_x}_{++}$, and vector $\lambda\in\mathbb{R}^{n_\phi} \geq 0$ such that
\begin{equation}\label{eq:model-stability}
    \tilde R_V^\top \begin{bmatrix}
        A_G^\top PA_G-P&A_G^\top PB_G\\
        B_G^\top PA_G&B_G^\top PB_G
    \end{bmatrix}
    \tilde R_V+\tilde R_\phi\begin{bmatrix}
        \Lambda&0\\0&-\Lambda
    \end{bmatrix}
    \tilde R_\phi\prec 0,
\end{equation}
where 
\begin{equation}\label{eq:model-stability-parameters}
    \tilde R_V=\begin{bmatrix}
        I_{n_x}&0\\
        \tilde N_{ux}&\tilde N_{uz}
    \end{bmatrix},
    \tilde R_\phi=\begin{bmatrix}
            \tilde N_{vx}&\tilde N_{vz}\\
            0&I_{n_\phi}
        \end{bmatrix}.
\end{equation}
Substituting \eqref{eq:model-stability-parameters} into \eqref{eq:model-stability}, and applying Schur complements, we obtain the following equivalent condition
\begin{equation}\label{eq:stability-nonconvex}
    \begin{bmatrix}
        P&0&A_G^\top +\tilde N_{ux}^\top B_G^\top & \tilde N_{vx}^\top \\
        0&\Lambda&\tilde N_{uz}^\top B_G^\top & \tilde N_{vz}^\top \\
        A_G+B_G\tilde N_{ux}&B_G\tilde N_{uz}&P^{-1}&0\\
        \tilde N_{vx}&\tilde N_{vz}&0&\Lambda^{-1}
    \end{bmatrix}\succ 0.
\end{equation}
Note here $A_G$ and $B_G$ are unknown. Introduce matrices $G_1\in\mathbb{R}^{T\times n_x}$, $G_2\in\mathbb{R}^{T\times n_\phi}$ satisfying
\begin{equation}\label{eq:data-parameterization}
    \begin{bmatrix}
        \tilde N_{ux}\\I_{n_x}
    \end{bmatrix}=
    \begin{bmatrix}
        U_0\\X_0
    \end{bmatrix}G_1,\begin{bmatrix}
        \tilde N_{uz}\\0
    \end{bmatrix}=
    \begin{bmatrix}
        U_0\\X_0
    \end{bmatrix}G_2.
\end{equation}
Using \citep[Theorem 1]{de2019formulas} and assuming Assumption \ref{ass:data-exciting} holds, we have
\begin{equation}\label{eq:data-representation}
    [B_G~A_G]\begin{bmatrix}
        U_0\\X_0
    \end{bmatrix}=X_1.
\end{equation}
Then we have
\begin{equation*}
     A_G+B_G\tilde N_{ux}=[B_G~A_G]\begin{bmatrix}
         U_0\\X_0
     \end{bmatrix}G_1\overset{}{=}X_1G_1,
\end{equation*}
and similarly, $B_G\tilde N{uz}=X_2G_2$. Substituting these into \eqref{eq:stability-nonconvex}, we obtain
\begin{equation}\label{eq:nonconvex}
    \begin{bmatrix}
        P&0&G_1^\top X_1^\top & \tilde N_{vx}^\top \\
        0&\Lambda&G_2^\top X_1^\top  & \tilde N_{vz}^\top \\
        X_1G_1&X_1G_2&P^{-1}&0\\
        \tilde N_{vx}&\tilde N_{vz}&0&\Lambda^{-1}
    \end{bmatrix}\succ 0.
\end{equation}
The constraint \eqref{eq:nonconvex} appears to be nonconvex. Multiplying $\text{diag}(P^{-1},\Lambda^{-1},I_{n_x},I_{n_\phi})$ on both sides of the matrix in \eqref{eq:nonconvex}, we obtain \eqref{eq:stability-LMI}, where $Q_1=P^{-1}$, $Q_2=\Lambda^{-1}$, $L_1=G_1Q_1$, $L_2=G_2Q_2$. To obtain \eqref{eq:data-representation-matrix}, multiply $Q_1$ on the right of both sides of $G_1$ and $G_2$ in \eqref{eq:data-parameterization}.
\end{proof}

\subsection{Safety Verification}
For safety verification, let $\mathcal{R}^k$ denote the exact reachable set, which is defined recursively as
\begin{equation}\label{eq:reach-set}
    \mathcal{R}^{k+1}=A_G\mathcal{R}^k+B_G\pi(\mathcal{R}^k),\quad k=0,\ldots,T-1,
\end{equation}
where $\mathcal{R}^0=\mathcal{X}\subseteq \mathcal{S}$. With a slight abuse of notation, $A_G\mathcal{R}^k:=\{A_Gx:x\in\mathcal{R}^k\}$, $\pi(\mathcal{R}^k):=\{\pi(x):x\in\mathcal{R}^k\}$. Then, the system is safe if and only if 
\begin{equation}\label{eq:safe-condition}
    \mathcal{R}^{k}\subseteq \mathcal{S}, \quad \forall k=1,\ldots,T.
\end{equation}
In practice, it is hard to compute the exact expression of $\mathcal{R}^k$, for $k=1,\ldots,T$, due to the nonlinearities in $\pi(x)$. Instead, one can compute (numerically) an outer approximation $\overline{\mathcal{R}}^k$, which satisfies $\mathcal{R}^k\subseteq \overline{\mathcal{R}}^k$ for $k=1,\ldots,T$ \citep{fazlyab2020safety}. Using this approximated reachable set, a sufficient condition for safety is 
\begin{equation}\label{eq:suff-safety-condition}
    \overline{\mathcal{R}}^k\subseteq \mathcal{S}, \quad \forall k=1,\ldots,T.
\end{equation}
Sufficiency can be seen by $\mathcal{R}^k\subseteq \overline{\mathcal{R}}^k\subseteq \mathcal{S},k=1,\ldots,T$. In order to efficiently compute $\overline{\mathcal{R}}^k$,  consider using a linear parameterization:
\begin{equation}
    \overline{\mathcal{R}}^k:=\left\{x\in\mathbb{R}^{n_x}:\bigcap_{i=1}^p{d_i}^\top (x-x^{safe})+\gamma_i^k\ge 0\right\},
\end{equation}
where for $i=1,\ldots,r$, $d_i\in\mathbb{R}^{n_x}$ is the same as the slopes used to define the safe set $\mathcal{S}$ in \eqref{eq:safe-set}, and $\gamma_i^k\in\mathbb{R}$ is a decision variable to be optimized. We set $\gamma_i^0=1$, for $i=1,\ldots,p$, such that $\overline{\mathcal{R}}^0=\mathcal{R}^0=\mathcal{X}$ as an initialization. Then, for $k=1,\ldots,T$, $\overline{\mathcal{R}}^k$ can be computed recursively by the following optimization problem
\begin{subequations}\label{eq:reach-opt-prob}
\begin{align}
    \min_{\{\gamma_i^k\}_{i=1}^p}~&\sum_{i=1}^p\gamma_i^k\\
    \mathrm{subject~to}~&\mathcal{H}_k:=\{A_Gw^0+B_G(W^{l+1}w^l+b^{l+1}):w^0\in\overline{\mathcal{R}}^{k-1},\eqref{eq:sector}\}\subseteq \overline{\mathcal{R}}^{k},
\end{align}
\end{subequations}
where $\overline{\mathcal{R}}^{k-1}$ is the solution from the last recursion, and $\gamma_i^k$ is the free decision variable from \eqref{eq:reach-set}. Here we have dropped the time argument $k$ for $w$ in the optimization problem, and substituted $x$ by $w^0$, and $\pi(x)$ by $W^{l+1}w^l+b^{l+1}$ from  \eqref{eq:nn}. In \eqref{eq:reach-opt-prob}, although the nonlinearities from the activation functions has been bounded by the offset sector using Lemma \ref{lem:sector}, the optimization problem is still hard to solve due to the set inclusion constraint. Both $\mathcal{H}^k$ and $\overline{\mathcal{R}}^k$ are \emph{semi-algebraic sets} as they are defined by real polynomials. Thanks to this property, \emph{Positivestellensatz} provides a tractable way to relax this set inclusion constraint into a sum-of-squares constraint, which is equivalent to a semi-definite constraint. To see this, the constraint $\mathcal{H}^k\subseteq \overline{\mathcal{R}}^k$ is equivalent to
\begin{equation}\label{eq:reach-empty-set-cond}
    \mathcal{H}^k\bigcap \{\overline{\mathcal{R}}^k\}^c=\emptyset,
\end{equation}
where $\{\overline{\mathcal{R}}^k\}^c$ is the closure of the complementary set of $\overline{\mathcal{R}}^k$, defined by
\begin{equation}\label{eq:compl-reach-set}
    \{\overline{\mathcal{R}}^k\}^c:=\left\{x\in\mathbb{R}^{n_x}:\bigcup_{i=1}^pd_i^\top (x-x^{safe})+\gamma_i\le 0\right\}.
\end{equation}
Then, \eqref{eq:reach-empty-set-cond} can be equivalently formulated as
\begin{equation}\label{eq:equiva-reach-cond}
    \mathcal{H}^k\bigcap\left\{x\in\mathbb{R}^{n_x}:d_i^\top (x-x^{safe})+1\le 0\right\}=\emptyset,\quad i=1,\ldots,p.
\end{equation}
For $k=2,\ldots,T$, the optimization problems \eqref{eq:reach-opt-prob} obtained using Positivestellensatz are as follows. 
\begin{subequations}\label{eq:psatz-reach-opt-model}
\begin{align}
\min~&\sum_{i=1}^p\gamma_i^k\\
\mathrm{subject~to}~&d_i^\top[A_Gw^0+B_G(W^{l+1}w^l+b^{l+1})-x^{safe}]+\gamma_i^k-\nonumber\\
&\sum_{i=1}^pg_i(w)[d_i^\top (w^0-x^{safe})+\gamma_i^{k-1}]\nonumber\\
&-\begin{bmatrix}
            v_\phi-v_*\\
            w_\phi-w_*
        \end{bmatrix}^\top
        \begin{bmatrix}
            -2A_\phi B_\phi\Lambda_i&(A_\phi+B_\phi)\Lambda_i\\
            (A_\phi+B_\phi)\Lambda_i&-2\Lambda_i
        \end{bmatrix}
        \begin{bmatrix}
            v_\phi-v_*\\
            w_\phi-w_*
        \end{bmatrix}\in\Sigma[w],\quad i=1,\ldots,p\label{eq:psatz-reach-constraint-model}\\
        &\gamma_i^k\in\mathbb{R},\quad g_i(w)\in\Sigma[w],\quad i=1,\ldots,p,\quad \lambda_i\in\mathbb{R}^{n_\phi}\ge 0,\quad \Lambda_i=\text{diag}(\lambda_i)
\end{align}
\end{subequations}
where we recall $w_\phi=[w^1,\ldots,w^l]^\top$, $A_\phi=\text{diag}(\alpha_\phi)$, $B_\phi=\text{diag}(\beta_\phi)$, and derive $v_\phi=[W^1w^0+b^1,\ldots,W^lw^{l-1}+b^l]^\top$ from \eqref{eq:nn}. The left hand side of constraint \eqref{eq:psatz-reach-constraint-model} is a polynomial in $w = [w^0,w^1,\ldots,w^l]^\top$. The decision variables represented by $g_i(w)$ are actually the coefficients of these polynomials. Given that the objective function is linear, \eqref{eq:psatz-reach-opt-model} is a semi-definite programming problem that can be solved efficiently by interior-point solvers.

The computation of $\overline{\mathcal{R}}^1$ is slightly different since $\overline{\mathcal{R}}^0=\mathcal{X}$, where $\mathcal{X}$ is the input set defined in \eqref{eq:input-set}. The computational program is similar to \eqref{eq:psatz-reach-opt-model}, the only modification is substituting $d_i$ by $c_i$, and $\gamma_i^{k-1}$ by one.

To solve the problem \eqref{eq:psatz-reach-opt-model} more efficiently, one can decompose it into $p$ problems, where every problem is to optimize $\gamma_i^k$, and subject to the $i^{th}$ constraint in \eqref{eq:psatz-reach-constraint-model}. This is because the $p$ constraints in \eqref{eq:psatz-reach-constraint-model} are independent to each other, and the objective function is in linear form, which is decomposable. 

Program \eqref{eq:psatz-reach-opt-model} is based on known $A_G$ and $B_G$. For the case that $G$ is unknown, $A_G$ and $B_G$ can be represented by data \eqref{eq:data-representation} into the program. The data-driven safety verification programs for $i=1,\ldots,p$ are given as follows.

\begin{subequations}\label{eq:psatz-reach-opt-data}
\begin{align}
\min~&\gamma_i^k\nonumber\\
\mathrm{subject~to}~&d_i^\top[X_1G_1 w^0+X_1G_2w^l+X_1G_3-x^{safe}]+\gamma_i^k-\nonumber\\
&\sum_{i=1}^pg_i(w)[d_i^\top (w^0-x^{safe})+\gamma_i^{k-1}]\nonumber\\
&-\begin{bmatrix}
            v_\phi-v_*\\
            w_\phi-w_*
        \end{bmatrix}^\top
        \begin{bmatrix}
            -2A_\phi B_\phi\Lambda&(A_\phi+B_\phi)\Lambda\\
            (A_\phi+B_\phi)\Lambda&-2\Lambda
        \end{bmatrix}
        \begin{bmatrix}
            v_\phi-v_*\\
            w_\phi-w_*
        \end{bmatrix}\in\Sigma[w],\label{eq:psatz-reach-constraint-data}\\
&\begin{bmatrix}
    0\\I_{nx}
\end{bmatrix}=\begin{bmatrix}
    U_0\\X_0
\end{bmatrix}G_1,\quad \begin{bmatrix}
    W^{l+1}\\0
\end{bmatrix}=\begin{bmatrix}
    U_0\\X_0
\end{bmatrix}G_2,\quad \begin{bmatrix}
    b^{l+1}\\0
\end{bmatrix}=\begin{bmatrix}
    U_0\\X_0
\end{bmatrix}G_3\label{eq:data-representation-safety}\\
\gamma_i^k\in\mathbb{R},& g_i(w)\in\Sigma[w], 0 \le \lambda\in\mathbb{R}^{n_\phi}, \Lambda=\text{diag}(\lambda),G_1\in\mathbb{R}^{K\times n_x}, G_2\in\mathbb{R}^{K\times n_\phi}, G_3\in\mathbb{R}^{K\times 1}
\end{align}
\end{subequations}
The program is convex as it only has sum-of-squares constraints and linear constraints, while the objective function is linear.
\begin{theorem}\label{th:safe}
    Consider system \eqref{eq:system} with an unknown plant $G$ and an NN controller $\pi$ as in \eqref{eq:nn}. Suppose $\mathcal{X}$ and $\mathcal{S}$ satisfy Assumption \ref{ass:safe}, data $X_0,X_1,U_0$ in \eqref{eq:data} satisfy Assumption \ref{ass:data-exciting}. Let $\alpha_\phi$, $\beta_{\phi}$, $v_*\in\mathbb{R}^{n_\phi}$ be given with $\alpha_\phi\le \beta_\phi$, and $w_*:=\phi(v_*)$. Assume $\phi$ element-wisely satisfies the offset sector $[\alpha_\phi,\beta_\phi]$, centered on the point $(v_*,w_*)$. Assume for $k=1,\ldots,T$, program \eqref{eq:psatz-reach-constraint-data} is feasible with the optimal solution denoted by $\{\gamma_i^k,g_i(w),\Lambda_i\}_{i=1}^p$. Then, System \eqref{eq:system} is safe in the time interval $[1,T]$ if
    \begin{equation}\label{eq:safety-condition}
        \gamma_i^k\le 1,\quad i=1,\ldots,p,\quad k=1,\ldots,T.
    \end{equation}
\end{theorem}

\begin{proof}
    From \eqref{eq:data-representation-safety} we obtain
    \[
    A_G=[B_G~A_G]\begin{bmatrix}
        0\\I_{nx}
    \end{bmatrix}=[B_G~A_G]\begin{bmatrix}
        U_0\\X_0
    \end{bmatrix}G_1=X_1G_1.
    \]
    Similarly, $B_GW^{l+1}=X_1G_2$, $B_Gb^{l+1}=X_1G_3$. Substituting $X_1G_1$ by $A_G$, $X_1G_2$ by $B_GW^{l+1}$, and $X_1G_3$ by $B_Gb^{l+1}$ into \eqref{eq:psatz-reach-constraint-data}, we deduce that \eqref{eq:psatz-reach-opt-data} is equivalent to \eqref{eq:psatz-reach-opt-model}. Therefore, $\overline{\mathcal{R}}^k$ for system \eqref{eq:system} can be computed by \eqref{eq:psatz-reach-opt-data} using data \eqref{eq:data}. We then prove that for $k=1,\ldots,T$, $\overline{\mathcal{R}}^k\subseteq \mathcal{S}$ provided by \eqref{eq:safety-condition}. \eqref{eq:safety-condition} indicates that 
    \[
    \{x\in\mathbb{R}^{n_x}:d_i^\top(x-x^{safe})+\gamma_i^k\ge 0\}\subseteq \{x\in\mathbb{R}^{n_x}:d_i^\top(x-x^{safe})+1\ge 0\}.
    \]
    Using \citep[Problem 2.33]{ashlock2020introduction} we have
    \[
        \overline{\mathcal{R}}^k=\bigcap_{i=1}^p\{x\in\mathbb{R}^{n_x}:d_i^\top(x-x^{safe})+\gamma_i^k\ge 0\}\subseteq \bigcap_{i=1}^p\{x\in\mathbb{R}^{n_x}:d_i^\top(x-x^{safe})+1\ge 0\}=\mathcal{S}.
    \]
This concludes the proof.
\end{proof}

For a given $T$, the number of programs to be solved for verifying safety is $pT$. Although the number grows linearly with $T$, verifying safety over a long time interval $[1,T]$ is challenging using the approximated reachability method. This is because the approximation error $|\mathcal{R}^k-\overline{\mathcal{R}}^k|$ accumulates and propagates with time. In the following section, we show results about verifying invariance for set $\mathcal{B}$. The results can be applied to verify safety for $T\to+\infty$, while not accumulating approximation errors.   

\subsection{Invariance Verification}
The invariance verification problem \ref{prob:invariance} can be regarded as a one-step reachability problem. To see this, let $\mathcal{R}^0=\mathcal{B}$. Then, $\mathcal{B}$ is invariant if $
\overline{\mathcal{R}}^1\subseteq \mathcal{B}$, where $\overline{\mathcal{R}}^1$ is solved from \eqref{eq:psatz-reach-opt-data} with $k=1$, $\mathcal{R}^0=\mathcal{B}$, $\mathcal{B}$ satisfies Assumption \ref{ass:invariance}. 

The program for invariance verification is proposed as follows.
\begin{subequations}\label{eq:psatz-inv-opt-data}
\begin{align}
\min~&\sum_{i=1}^p\gamma_i\\
\mathrm{subject~to}~&e_i^\top[X_1G_1 w^0+X_1G_2w^l+X_1G_3-x^{safe}]+\gamma_i-\nonumber\\
&g_i(w)[e_i^\top (w^0-x^{safe})+1]\nonumber\\
&-\begin{bmatrix}
            v_\phi-v_*\\
            w_\phi-w_*
        \end{bmatrix}^\top
        \begin{bmatrix}
            -2A_\phi B_\phi\Lambda_i&(A_\phi+B_\phi)\Lambda_i\\
            (A_\phi+B_\phi)\Lambda_i&-2\Lambda_i
        \end{bmatrix}
        \begin{bmatrix}
            v_\phi-v_*\\
            w_\phi-w_*
        \end{bmatrix}\in\Sigma[w],i=1,\ldots,p\label{eq:psatz-inv-constraint-data}\\
        &\begin{bmatrix}
    0\\I_{nx}
\end{bmatrix}=\begin{bmatrix}
    U_0\\X_0
\end{bmatrix}G_1,\begin{bmatrix}
    W^{l+1}\\0
\end{bmatrix}=\begin{bmatrix}
    U_0\\X_0
\end{bmatrix}G_2,\begin{bmatrix}
    b^{l+1}\\0
\end{bmatrix}=\begin{bmatrix}
    U_0\\X_0
\end{bmatrix}G_3\\
        \gamma_i^k\in\mathbb{R},&g_i(w)\in\Sigma[w],\lambda\in\mathbb{R}^{n_\phi}\ge 0,\Lambda=diag(\lambda), G_1\in\mathbb{R}^{K\times n_x},G_2\in\mathbb{R}^{K\times n_\phi},G_3\in\mathbb{R}^{K\times 1}
\end{align}
\end{subequations}

\begin{theorem}\label{th:invariance}
    Consider system \eqref{eq:system} with an unknown plant $G$ and an NN controller $\pi$ as in \eqref{eq:nn}. Suppose $\mathcal{B}$ satisfies Assumption \ref{ass:invariance}, data $X_0,X_1,U_0$ in \eqref{eq:data} satisfy Assumption \ref{ass:data-exciting}. Let $\alpha_\phi$, $\beta_{\phi}$, $v_*\in\mathbb{R}^{n_\phi}$ be given with $\alpha_\phi\le \beta_\phi$, and $w_*:=\phi(v_*)$. Assume $\phi$ element-wisely satisfies the offset sector $[\alpha_\phi,\beta_\phi]$, centered on the point $(v_*,w_*)$. Assume program \eqref{eq:psatz-inv-opt-data} is feasible, and the optimal solution is denoted by $\{\gamma_i,g_i(w),\Lambda_i\}_{i=1}^q$. Then, the set $\mathcal{B}$ is invariant under vector field $A_Gx+B_G\pi$ if
    \begin{equation}
        \gamma_i\le 1,i=1,\ldots,q.
    \end{equation}
\end{theorem}
\begin{proof}
    The proof is similar to that for Theorem \ref{th:safe}.
\end{proof}

From \citep[Lemma 1]{wang2023safety}, if $\mathcal{X}\subseteq \mathcal{B}\subseteq 
\mathcal{S}$, and $\mathcal{B}$ is invariant, system \eqref{eq:system} is safe.

\section{Simulation Results}
\label{sec:simulation}
\subsection{Vehicle Cruise Control}
% controller trained by imitation learning

Consider the vehicle lateral dynamics problem taken from \citep[Section IV. B]{yin2021stability}, which is a forth-order linear system. The state of the system is $x=[e,\dot e,e_\theta,\dot{e}_\theta]^\top$, where $e$ is the perpendicular distance to the lane centre, and $e_\theta$ is the angle between the tangent to the straight section of the road and the projection of the vehicle's longitudinal axis. The disturbance $c$ is set to zero. 

In that paper, an NN controller $\pi$ is trained under certain parameters of the dynamics. However, these parameters can be uncertain, or even time-varying in practice. Data-driven verification can be highly efficient for this case. In our simulation, we consider sampling input-state data for the plant, where all the parameters, such as the longitudinal velocity, front cornering stiffness, etc, are subject to unknown Gaussian noise. The solution of data-driven method is compared with that of the model-based method. We also use different sample size of data and observe how it affects the result. 

% The data are collected with $K=100$.
\subsubsection{Stability Verification}
From Theorem \ref{th:stability}, we know the feasibility of \eqref{eq:data-stability-matrix}, where $Q_1$ is a positive definite matrix and $Q_2$ is a diagonal matrix, leads to the stability of system. To find a region of attraction (ROA) as large as possible for this vehicle system, the stability verification is transformed into an semidefinite programming (SDP) problem that minimizes the trace of $Q_1$ with \eqref{eq:data-stability-matrix} as constraints. 

For the data-driven method, we notice that the problem is infeasible when the sample size $K<5$. This is because the rank condition \eqref{eq:exciting} is not satisfied: the dimension of system state $n_x$ and dimension of system input $n_u$ are 4 and 1 respectively. Therefore, the smallest $K$ that satisfies rank condition is 5. When $K\geq 5$, the problem is feasible and the solutions are almost the same. For the model-based method, the problem is always feasible and the obtained solution is the same as that obtained from data-driven method.

\begin{figure}
\centering
  \includegraphics[width=0.8\linewidth]{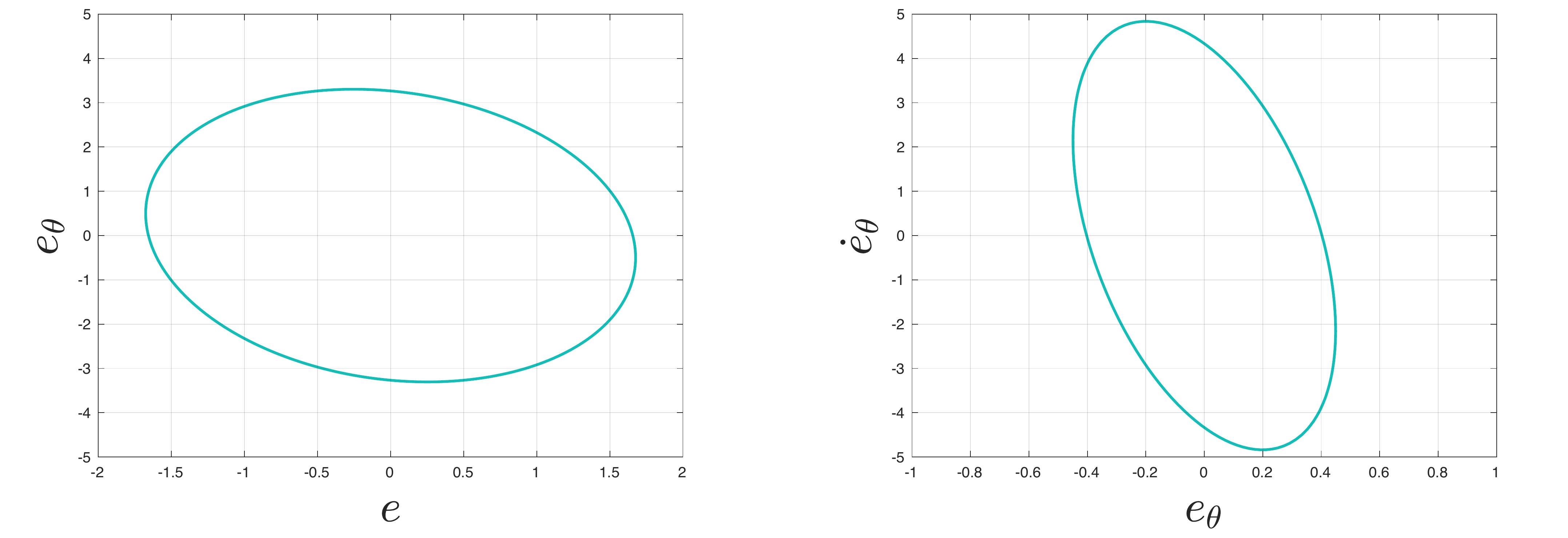}
  \caption{The ROAs of NN controller.}
  \label{fig:ROA}
\end{figure}

The region of attraction (ROA) is defined by \citep{yin2021imitation} as:
\begin{equation}
    \mathcal{E}(P) := \{x \in \mathbb{R}^{n_x} : x^TPx \leq 1\}
\end{equation}
Combining this definition with the solution for $Q_1$, which is the inverse of $P$, we can get the ROAs of the system as Fig.\ref{fig:ROA} shows. As the dimension of state $x$ is four, we set $[e_\theta,\dot{e}_\theta]$ and $[e,\dot{e}]$ as $[0,0]$ respectively and project their ROAs onto the corresponding planes. It is obvious that they are contained by the range of the states. 

\subsubsection{Safety Verification}
We consider the input set $\mathcal{X}$ defined by inequalities $-0.1\le e \le 0.1$, $-0.1\le \dot e\le 0.1$, $-\pi/60\le e_\theta \le \pi/60$, $-\pi/60\le \dot e_\theta \le \pi/60$. The safe set $\mathcal{S}$ is defined by inequalities $-2\le e\le 2$. By solving the data-driven safety verification program \eqref{eq:psatz-reach-opt-data}, we obtain $\gamma_1^1=0.051$, $\gamma_2^1=0.051$, $\gamma_1^2=0.4273$, $\gamma_2^2=0.566$, $\gamma_1^3=0.9254$, $\gamma_2^3=1.2199$. Using Theorem \ref{th:safe}, this demonstrates that the system \eqref{eq:system} is safe for $k=1,2$. 

\subsection{Pendulum}\label{sec:pendulum}
Consider an inverted pendulum with dynamics
\begin{equation}\label{eq:pendulum}
    \dot x=\begin{bmatrix}
        \omega\\
        \frac{3g}{2l}\sin\theta-k\omega
    \end{bmatrix}
    +\begin{bmatrix}
        0\\\frac{3}{ml^2}u
    \end{bmatrix},
\end{equation}
where $x=[\theta~\omega]\in\mathbb{R}^2$ is the state, $u\in\mathbb{R}$ is the torque which controls the pendulum. The parameters $k=0.8$ represents the damping ratio for air resistance, $g=10.0$ is the gravity constant, $l=1.0$ is the length of the pendulum and $m=1.0$ is the weight. Controller $u$ is trained to stabilize the system around the equilibrium $x_*=(3.14,0)$ by the Soft Actor-Critic (SAC) method. Specifically, $u$ is a neural network with one hidden layer of 128 neurons, and tanh is used as the activation functions. Even though such a system is nonlinear, our proposed data-driven stability verification conditions \eqref{eq:data-stability-verification} can still be applied since the nonlinear dynamics can be linearized around the equilibrium $x_*$. One benefit of our method is that explicit linearization is not necessary, as we only need to sample data around the equilibrium $x_*$. 

In the experiment, the system dynamics \eqref{eq:pendulum} is digitalized with a time interval $dt=0.02$. Length of trajectory $K=5$, the open loop control is randomly sampled by $u(t)\in[-0.01,0.01]$, for $t=1,\ldots,5$ and the initial state $x(0)$ is randomly sampled by $\theta(0)\in[3.13,3.15],\omega(0)\in[-0.01,0.01]$. By solving an SDP using our proposed data-driven stability verification conditions \eqref{eq:stability-LMI}, we have $Q_1=\begin{bmatrix}
    10.338&-1.075\\-1.075&0.618
\end{bmatrix}\succ 0$. The RL controlled pendulum is verified to be locally stable around $x_*$ using the local Lyapunov function $x^\top Q_1^{-1}x$.

\section{Conclusion}
\label{sec:conclusion}
In this paper, we propose model-free stability and safety verification programs for neural network controlled systems. Our method does not require one to identify the system dynamics from input-state data, but directly uses these data into the verification program. In the future, we will explore extending the method to noise corrupted data for robust verification.

% Acknowledgments---Will not appear in anonymized version
\acks{}

\bibliography{ref}

\end{document}